\documentclass[10pt]{amsart}

\newcommand{\ra}{\rightarrow}
\newcommand{\R}{\mathbb{R}}
\newcommand{\T}{\mathbb{T}}
\newcommand{\Z}{\mathbb{Z}}
\newcommand{\N}{\mathbb{N}}

\newcommand{\ep}{\epsilon}

\newcommand{\8}{\infty}
\newcommand{\beq}{\begin{equation}}
\newcommand{\eeq}{\end{equation}}
\newcommand{\be}{\begin{eqnarray}}
\newcommand{\ee}{\end{eqnarray}}
\newcommand{\nn}{\nonumber}
\newcommand{\dfr}{\mbox{{\rm Diff}$^r(M)$}}
\newcommand{\Vol}{\mbox{\rm Vol}}
\newcommand{\PH}{\mbox{\rm PH}}
\newtheorem{thm}{Theorem}[section]
\newtheorem{lem}[thm]{Lemma}
\newtheorem{prop}[thm]{Proposition}
\newtheorem{dfn}[thm]{Definition}

\usepackage{color}

\def\diam{\mbox{\rm diam}}

\begin{document}

\title[entropy]{Topological Entropy and Partially Hyperbolic Diffeomorphisms}
\author{Yongxia Hua}
\address{Department of Mathematics, Northwestern University,
Evanston, Illinois 60208}
\email{hua@math.northwestern.edu}
\author{Radu Saghin}
\address{Department of Mathematics, University of Toronto, Toronto,
  Ontario, Canada M5S 2E4}
\email{rsaghin@fields.utoronto.ca}
\author{Zhihong Xia}
\address{Department of Mathematics \\ Northwestern University \\
Evanston, Illinois 60208}
\email{xia@math.northwestern.edu}
\thanks{Research supported in part by National Science Foundation.}
\date{July 24, 2006, draft version}

\begin{abstract}
  We consider partially hyperbolic diffeomorphisms on compact
  manifolds where the unstable and stable foliations stably carry some
  unique non-trivial homologies. We prove the following two results:
  if the center foliation is one dimensional, then the topological
  entropy is locally a constant; and if the center foliation is two
  dimensional, then the topological entropy is continuous on the set
  of all $C^\8$ diffeomorphisms. The proof uses a topological
  invariant we introduced; Yomdin's theorem on upper semi-continuity;
  Katok's theorem on lower semi-continuity for two dimensional systems
  and a refined Pesin-Ruelle inequality we proved for partially
  hyperbolic diffeomorphisms.
\end{abstract}

\maketitle

\section{Introduction and Main Results}

One of the fundamental invariants in topological dynamics is
topological entropy. However, the entropy is a very hard to compute
and its dependence on the map itself and it continuity properties are
very delicate.  For uniformly hyperbolic diffeomorphisms on a compact
manifold, the topological entropy is locally a constant. i.e., it
remains the same under small perturbations. This is due to the
structural stability of uniformly hyperbolic systems. In other words,
the entropy is stable for hyperbolic systems. The first question we
ask in this paper is the following: besides the structurally stable
systems, are there any other system where the topological entropy are
stable under perturbations? We will show that the answer to this
question is yes and there are classes of partially hyperbolic systems
with one dimensional centers where the topological entropy is locally
constant.

Our next question is: when is the topological entropy continuous?
This is a difficult problem. In general, entropy is not continuous
for a $C^1$ diffeomorphism. Yomdin \cite{Yo} proved that it is
upper semi-continuous for the class of $C^{\infty} (M)$
diffeomorphisms for any compact manifold.  For $\dim(M) = 2$,
Katok \cite{Po} showed that the the entropy is lower
semi-continuous for the $C^{1+{\alpha}}$, $\alpha >0$
diffeomorphisms on $M$. Combining these two, we have the
continuity of topological entropy for $C^\8$ diffeomorphisms on
compact surfaces.

In this paper we will show that for a large class of partially
hyperbolic $C^\8$ diffeomorphisms with two dimensional center
foliations the topological entropy is continuous. Besides requiring
that the center foliation has dimension two or less, we also require
that the stable and unstable foliations carry certain homological
information of the manifold. We will give a detailed definition later.

Let $M$ be a compact Riemannian manifold and let $f \in \PH^\8(M)$ be
the set of $C^\8$ partially hyperbolic diffeomorphisms on $M$. $f$ is
said to be {\it partially hyperbolic} if for every $x\in M$ the
tangent space at $x$ admits an invariant splitting
\[
T_xM=E^s(x)\oplus E^c(x)\oplus E^u(x)
\]
into {\it strongly stable} $E^s(x)=E^s_f(x)$, {\it central}
$E^c(x)=E^c_f(x)$, and {\it strongly unstable} $E^u(x)=E^u_f(x)$
subspaces and  there exist numbers $c_1 >1$ and
\[
0<\lambda_s < \lambda'_c \leq 1 \leq \lambda''_c < \lambda_u
\]
such that for every $x\in M$ and all $i \in \N$,
\be
v\in E^s(x)&\Rightarrow& \|d_xf^i(v)\| \leq \, c_1\lambda_s^i \|v\|, \nn \\
v\in E^c(x)&\Rightarrow& c_1^{-1}(\lambda'_c)^i\,\|v\| \leq \|d_xf^i(v)\|
                   \leq  c_1 (\lambda''_c)^i \, \|v\|, \label{eqn10}\\
v\in E^u(x)&\Rightarrow& \,c_1^{-1} \lambda_u^i\|v\| \leq
\|d_xf^i(v)\|. \nn
\ee
We denote the set of all $C^r$ partially hyperbolic diffeomorphisms by
$\PH^r(M)$.

\vspace{1ex}
We can state our main results of the paper:

\begin{thm}
  Let $M$ be a compact Riemannian manifold and let $f \in \PH^\8(M)$
  be the set of $C^\8$ partially hyperbolic diffeomorphisms on
  $M$. Assume that (1).\ the dimension of the center foliation is
  $2$
  or less; (2).\ the strong stable and strong unstable foliation
  stably carry some unique non-trivial homologies, then the
  topological entropy $h_{top}: \PH^\8(M) \ra \R$ is continuous at
  $f$.

  Furthermore, if the center foliation has dimension $1$, then
  $h_{top}$ is a constant in a small neighborhood of $f \in
  \PH^{1}(M)$.
\end{thm}

We will define the homologies carried by a foliation. If the stable
manifold and unstable manifold are one dimensional, then the
homological condition is a condition on the homotopy class of the
map. For higher dimensional case, we believe that the same is true,
but we are not able to prove that. We do have some open conditions
that one can verify. The theorem is not true in general without the
assumptions on the homology. In the last section of this paper, we
will give some examples where the theorem fails without such
assumptions.

Even though our main results are about the topological entropy, its
proof, however, relies on smooth ergodic theory, Lyapunov exponents
and measure theoretic entropy. In particular, we proved a
refined Pesin-Ruelle inequality for partially hyperbolic
diffeomorphisms. This is stated in Theorem \ref{thmpr}, which is of
interest in its own right.

\section{Currents, topological and geometric growth}

In this section, we will define some topological invariants for
diffeomorphisms with uniformly expanding (or contracting)
foliations. The topological invariant was introduced in Saghin and Xia
\cite{SX06}, where one can found more details and some other
applications of the invariant.  We will relate the volume growth of an
expanding invariant foliation with this topological invariant. If an
invariant foliation carries certain nontrivial homological information
of the manifold, which we will make precise later, then the volume
growth, which is harder to track, is exactly the same as the
topological growth. The topological growth can be easily calculated by
actions induced by the map on the homology of the manifold.

Let $M$ be an $n$-dimensional compact Riemannian manifold. Let $f \in
\dfr$ be diffeomorphism on $M$. Let $W$ be a $k$-dimensional foliation
of $M$, invariant under $f$, i.e., $f$ maps leaves of $W$ to
leaves. This invariant foliation will be, for the purpose of this
paper, the strongly stable and strongly unstable foliations. We first
define volume growth of $f$ on leaves of $W$. For any $x \in M$, Let
$W(x)$ be the leave through $x$ and let $W_r(x)$ be the $k$
dimensional disk on $W(x)$ centered at $x$, with radius $r$.

Let
$$\chi_W(x, r) = \lim \sup_{n \ra \8} \frac{1}{n} \ln (\Vol (f^n(W_r(x))))$$
$\chi_W(x, r)$ is the volume growth rate of the foliation at $x$. Let
$$\chi_W(f)  = \sup_{x \in M} \chi_W(x, r)$$
Then, $\chi_W(f)$ is the maximum volume growth rate of $W$ under $f$.
Obviously, the quantity $\chi_W(f)$ is independent of $r$.

We search for conditions such that $\chi_W(x, r)$ is independent of
both $x$ and $r$. This is not true in general, there need to be
certain topological conditions for this to hold. If $W$ is
exponentially expanding under the iterates of $f$, we can formulate
this condition in terms of the action $f$ induces on the homology of
$M$. To simply put it, if $f$ is partially hyperbolic with $W$ being
part of unstable manifolds, then we requires that $f_*: H_*(M, \R) \ra
H_*(M, \R)$ be partially hyperbolic in a compatible way.

More precisely, let $W$ be an $f$-invariant $k$-dimensional
foliation. We assume that $W$ is orientable and we will fix an
orientation for $W$. Furthermore, we assume that the leaves of $W$
have exponential growth under $f$. i.e., there are constants $\lambda
>1$ and $c_2 >0$ such that
$$|df_x^n v| \geq c_2 \lambda^n |v|$$ for all $x \in M$, all $v \in
T_xW(x)$ and all $n \in \N$, where $W(x)$ is the leaf of $W$
through the point $x$. Let $W_r(x)$ be the ball of radius $r$
centered at $x$ on the leaf $W(x)$. For any positive integer, we
define the {\it currents}:
\be C_n(\omega) = \frac{1}{\Vol(f^n(W_r(x)))}
\int_{f^n(W_r(x))}, \label{eqn5} \ee for any $k$-form $\omega$ on $M$. These
currents depend on $x$ and $r$. The currents are uniformly bounded so
there must be subsequences with weak limits. Let $C$ be such a limit,
i.e., we have a sequence $n_i \ra \8$ such that for any $k$-form
$\omega$ we have $\lim_{i \ra \8} C_{n_i}(\omega) = C(\omega)$.

A current $C$ is said to be {\it closed}\/ if for any exact $k$-form
$\omega = d \alpha$, we have $C(\omega) = C(d \alpha) =0$. If $C$ is
closed, it has a homology class $[C]= h_C \in H_k(M, \R)$. This homology
class is nontrivial if there exist a closed $k$-form $\omega$ such
that $C(\omega) \neq 0$.

We would like to investigate the conditions under which the
subsequential limits of the currents $C_n$ is closed. In general,
$C_n$ itself is not closed. We believe that it can be approximated by
a closed one for large $n$. From Stokes' Theorem, we have: \be
C_n(\omega) &=& \frac{1}{\Vol(f^n(W_r(x)))}
\int_{f^n(W_r(x))} d \alpha \nn \\
&=& \frac{1}{\Vol(f^n(W_r(x)))} \int_{W_r(x)}(f^*)^n d \alpha \nn \\
&=& \frac{1}{\Vol(f^n(W_r(x)))} \int_{ \partial W_r(x)} (f^*)^n \alpha
\nn \ee

There are reasons to believe that the above sequence always approaches
$0$ as $n \ra \8$. i.e., every subsequential limit of the currents
$C_n$ is closed. Nevertheless, we are not able to show this as of now.
However, we can indeed show that this is true in many cases.

The first case is when the dimension of the foliation is one. In
this case, $\alpha$ is a real valued function and hence $\int_{
\partial W_r(x)} (f^*)^n \alpha$ is the difference of that
function evaluated at the two end points of $f^n(W_r(x))$ and
therefore it is uniformly bounded. Thus $C_n(\omega) \ra 0$ as $n
\ra \8$.

Another case is that when $f$ is close to a linear map on the torus
$\T^n$ and $W$ is any of the expanding foliations close to the linear
one. We will consider this case in more details later.

In general, we have the following simple proposition, whose proof is
straightforward.

\begin{prop}
  Let $J_k(x)$ be the Jacobian of $f$ restricted to the unstable
  subspace of $x$ and let $J_{k-1}(x)$ be the maximal Jacobian on
  $k-1$-dimensional subspace at $x$. If $J_k(x) > J_{k-1}(x)$ for all
  $x \in M$, then all subsequential limits of $\{C_n\}$ in equation
  \eqref{eqn5} are closed.
\end{prop}

The Jacobian condition in the above proposition is an open condition.

\begin{dfn}
  We say that a $k$-dimensional invariant foliation $W$ carries a
  non-trivial homology $h_C \in H_k(M, \R)$ if the currents $C_n$
  defined above has a closed subsequential limit $C$ and $h_C = [C]
  \neq 0$.

  We say that a $k$-dimensional invariant foliation $W$ carries a
  unique non-trivial homology (up to rescale) if all subsequential
  limits of the currents $C_i$ are closed and the homologies it carries
  are unique up to scalar multiplication and are uniformly bounded
  away from zero, for all $x \in M$ and all $r > 0$.
\end{dfn}

A current is non-trivial if there is a closed $k$-form $\omega$ such
that $C(\omega) \neq 0$. The homology class of a non-trivial closed
current is non-trivial. One way to show that the closed current $C$ is
non-trivial is to show that there is a closed $k$-form $\omega$ such
that $\omega$ is non-degenerate on $T_xW(x)$ for any $x \in M$. This
condition implies that the integral of $\omega$ over any oriented
segment of $W$ is nonzero. i.e.,
$$\int_D \omega \neq 0$$ for any piece $D$ on a leaf of the foliation $W$,
with its orientation inherited from the leaf. We may assume that the
integral is positive by choosing $-\omega$ if necessary. When we have
a non-degenerate $k$-form on the leaves of $W$, by compactness of the
manifold, there exists a constant $c_2 > 1$ such that
$$ c_2^{-1} \Vol(D) \leq \int_{D}\omega \leq c_2
\Vol(D) $$ for any segment $D$ on the leaves of $W$ and therefore
$$ c_2^{-1} \Vol(f^n(W_r(x))) \leq \int_{f^n(W_r(x))}\omega \leq c_2
\Vol(f^n(W_r(x))) $$ This implies that $C(\omega) >0$.

Assume that an invariant foliation $W$ carries a non-trivial homology
and let $h_C=[C] \in H_k(M, \R)$, where $C$ is the current as defined
above. The next proposition shows that $h_C$ is actually an eigenvector
of the induced linear map by $f$ on the homology of $M$.

\begin{prop} \label{prop10}
  Let $W$ be a $k$-dimensional invariant foliation that carries a unique
  non-trivial homology $h_C$. Then $h_C$ is an eigenvector of the induced
  linear map: $$f_*: H_k(M, \R) \ra H_k(M, \R).$$
\end{prop}

\begin{proof}
First we observe that the map $f$ naturally induces an action on the
currents, defined by: $$f_*C(\omega) = C(f^*\omega)$$ for any $k$
current $C$ and $k$ form $\omega$. Obviously, if $C$ is closed, then
$f_*C$ is closed too and $$[f_*C] = f_* h_C \in H_k(M. \R).$$

Let current $C$ be a subsequential limit of $C_n(x, r)$, then
$$C(\omega) = \lim_{i \ra \8} \frac{1}{\Vol(f^{n_i}(W_r(x)))}
\int_{f^{n_i}(W_r(x))}\omega,$$ for any $k$-form on $M$. Therefore \be
(f_*C)(\omega) &=& \lim_{i \ra \8} \frac{1}{\Vol(f^{n_i}(W_r(x)))}
\int_{f^{n_i}(W_r(x))}f^*\omega \nn \\
&=& \lim_{i \ra \8} \frac{1}{\Vol(f^{n_i}(W_r(x)))}
\int_{f^{(n_i+1)}(W_r(x))} \omega \nn  \\
&=& \lim_{i \ra \8}
\frac{\Vol(f^{(n_i+1)}(W_r(x)))}{\Vol(f^{n_i}(W_r(x)))} \cdot
\frac{1}{\Vol(f^{(n_i+1)}(W_r(x)))} \int_{f^{(n_i+1)}(W_r(x))} \omega
\nn \ee Since the ratio
${\Vol(f^{(n_i+1)}(W_r(x)))}/{\Vol(f^{n_i}(W_r(x)))}$ is uniformly
bounded, both from above and away from zero, there is a convergent
subsequence. Without loss of generality, we may assume that the
sequence actually converges and there is a constant $\lambda >0$ such
that $$ \lim_{i \ra
  \8}\frac{\Vol(f^{(n_i+1)}(W_r(x)))}{\Vol(f^{n_i}(W_r(x)))} = \lambda$$
This implies that $f_*C/\lambda$ is also a subsequential limit of the
current $C_n(x, r)$. Since $W$ carries a unique non-trivial homology,
this limit must be a scalar multiple of $C$. Therefore, there is a
constant $c_3$ such that
we have $f_*C \lambda^{-1} = c_3  C$. This implies that
$$f_*h_C = c \lambda h_C$$
i.e., $h_C$ is an eigenvector of $$f_*: H_k(M, \R) \ra H_k(M, \R)$$
with corresponding eigenvalue $c_3 \lambda$.

This proves the proposition.
\end{proof}

Let $\lambda_W$ be the eigenvalue of $f_*$ corresponding to the
eigenvector $h_C$, as in the above proposition. We call $\lambda_W$
the {\it topological growth} of the foliation $W$. We will see below
that the topological growth and the volume growth are the same for a
foliation that carries a unique non-trivial homology, except that the
volume growth we defined here is an exponent, while the topological
growth is a multiplier.

\begin{prop}\label{1}
  Let $W$ be a hyperbolic invariant foliation that carries a unique
  non-trivial homology $h_W$. Let $\lambda_W$ be the topological growth
  of the foliation. Then the volume growth defined before,
  \begin{equation*} \chi_W(f) = \ln \lambda_W \end{equation*} for any
  $x \in M$ and any $r >0$.
\end{prop}

\begin{proof}
  The volume of a piece of leaf in a foliation depends on the
  Riemannian metric defined on $M$. So in general, the volume does not
  grow uniformly with each iteration. We will need to rescale the
  volume at each step so that there will be uniform growth. Let $h_W
  \in H_k(M, \R)$ be a homology carried by $W$. For any $x \in M$ and
  $r >0$, we choose a sequence of numbers $d_i$, $i \in \N$ such
  that $$\lim_{i \ra \8} d_i C_i = C \; \mbox{ and } \; [C] =
  h_W/||h_W||.$$ This is possible by the uniqueness of homologies
  carried by the foliation. Moreover, there are numbers $0 <c_4 \leq
  c_5$ such that $d_i$ can be chosen with $c_5^{-1} \leq d_i \leq
  c_4^{-1}$. Therefore, \be (f_*C)(\omega) &=& \lim_{i \ra \8}
  \frac{d_i}{\Vol(f^{i}(W_r(x)))}
  \int_{f^{i}(W_r(x))}f^*\omega \nn \\
  &=& \lim_{i \ra \8} \frac{d_i}{\Vol(f^{i}(W_r(x)))}
  \int_{f^{(i+1)}(W_r(x))} \omega \nn  \\
  &=& \lim_{i \ra \8}
  \frac{\Vol(f^{(i+1)}(W_r(x)))/d_{i+1}}{\Vol(f^{i}(W_r(x)))/d_i}
  \cdot \frac{d_{i+1}}{\Vol(f^{(i+1)}(W_r(x)))}
  \int_{f^{(i+1)}(W_r(x))} \omega
  \nn  \\
  &=& \lim_{i \ra \8}
  \frac{\Vol(f^{(i+1)}(W_r(x)))/d_{i+1}}{\Vol(f^{i}(W_r(x)))/d_i}
  \cdot C(\omega) \nn \ee Therefore $$\lim_{i \ra \8}
  \frac{\Vol(f^{(i+1)}(W_r(x)))/d_{i+1}}{\Vol(f^{i}(W_r(x)))/d_i} =
  f_*C(\omega)/C(\omega) = \lambda_W$$

This implies that \be \chi_W(x,
  r) &=& \lim \sup_{i \ra \8} \frac{1}{i} \ln
  (\Vol (f^i(W_r(x)))) \nn \\
&=& \lim \sup_{i \ra \8} \frac{1}{i} \ln
  \Vol(f^i(W_r(x))) \nn \\
&=& \lim \sup_{i \ra \8} \frac{1}{i} \ln
(d_i^{-1}\Vol(f^i(W_r(x)))) \nn \\
  &=& \lim \sup_{i \ra \8} \frac{1}{i} \ln \left(d_0^{-1} \Vol(W_r(x))
    \cdot (\prod_{j=1}^i
  \frac{d_i^{-1}\Vol(f^{i}(W_r(x)))}{d_{i-1}^{-1}\Vol(f^{(i-1)}(W_r(x)))})
\right)
   \nn \\
&=& \lim \sup_{i \ra \8} \frac{1}{i} \sum_{j=1}^i \left(\ln
  \frac{d_i^{-1}\Vol(f^{i}(W_r(x)))}{d_{i-1}^{-1}\Vol(f^{(i-1)}(W_r(x)))}
\right) \nn \\
&=& \ln \lambda_W \nn \ee
Here we used the elementary fact that if $\lim_{i \ra \8} a_i =a$,
then $$\lim_{i \ra \8} \frac{1}{i} \sum_{j=1}^i a_i = a.$$

This proves the proposition.
\end{proof}

The next proposition discusses the situation where a foliation carries
more than one non-trivial homologies.

\begin{prop} \label{prop30}
  Let $W$ be a hyperbolic invariant foliation and let $H \subset
  H_k(M, \R)$ be the set
  of non-trivial homologies carried by $W$. Then $H$ spans a linear
  space, invariant under $$f_*: H_k(M, \R) \ra H_k(M, \R).$$
\end{prop}

\begin{proof}
We first observe that $H \subset H_k(M, \R)$ is a bounded and closed
set. Let $h \in H$ be a homology carried by the foliation $W$. It follows
from the proof of Proposition \ref{prop10} that there exists a
constant $c_6 >0$ such that $f_*h/c_6$ is also carried by $W$. The
proposition follows.
\end{proof}

Suppose that the foliation $W$ is one dimensional, then every
subsequential limit of the currents is closed, as we have shown. In
this case every limit defines a homology class. If we furthermore
assume that there is a non-degenerate closed 1-form $\omega$ on the
leaves of the foliation, then $W$ carries a nontrivial homology. Then
this homology class is non-trivial.

Another class of maps that we would like to consider are the ones on
$n$-torus $\T^n$ close to a linear map. Consider an $n \times n$
matrix $A$ with integer entries and with determinant one. The matrix
$A$ induces a linear toral automorphism: $T_A: \T^n = \R^n / \Z^n \ra
\T^n$ defined by $T_Ax = Ax \mod \Z^n$. If all eigenvalues are away
from the unit circle, then $T_A$ is a hyperbolic toral
automorphism. If the eigenvalues of $A$ are mixed, with some on the
unit circle and some away from unit circle, then $T_A$ is partially
hyperbolic.

In both hyperbolic and partially hyperbolic cases, let $E^u$ be
the unstable distribution of the $T_A$ on $\T^n$. At each point $x
\in \T^n$, $E^u(x) \subset T_x\T^n$ is the unstable subspace for
$dT_A: T_x\T^n \ra T_x\T^n$. Let $W^u$ be the unstable foliation
generated by $E^u$, $W^u$ is a hyperplane in $\T^n$. It is easy to
see that the currents $C_n$ converges to a unique closed current
$C$ and $C$ is non-trivial. Moreover, the eigenvalue corresponding
to $h_C$ for the map $f_*: H_*(\T^n, \R) \ra H_*(\T^n, \R)$ is the
product of all eigenvalues outside of the unit circle. i.e.,
$\lambda_W = \prod_{|\lambda_i| >1} \lambda_i$.

Let $f$ be a map close to $T_A$. We claim that all the subsequential
limits of the currents $C_n$ are closed. This is because that the
Jacobian for $f$ on the $k$-dimensional is close to the Jacobian for
$T_A$, which is equal to $(\prod_{|\lambda_i| >1}
\lambda_i)^n$. Therefore the $k$ dimensional volume
$\Vol(f^n(W_r(x)))$ grows with a factor close to $(\prod_{|\lambda_i|
  >1} \lambda_i)^n$. However, any $(k-1)$-form $(f^*)^n\alpha$ grows
approximately at the rate of the product of $k-1$
eigenvalues. Therefore,
\be
C_n(\omega) &=& \frac{1}{\Vol(f^n(W_r(x)))}
\int_{f^n(W_r(x))} d \alpha \nn \\
&=& \frac{1}{\Vol(f^n(W_r(x)))} \int_{ \partial W_r(x)} (f^*)^n \alpha
\ra  0  \nn \ee
as $ n \ra \8$.

It is also easy to see that every subsequential limit of the currents
is non-trivial. For the linear map, there is a coordinate plane
$(x{n_1}, x_{n_2}, \ldots, x_{n_k})$ such that the orthogonal
projection of the unstable space $E^u$ to the plane is
nondegenerate. Then the $k$-form $\omega^k = dx_{n_1}\wedge \cdots
\wedge dx_{n_k}$ is nondegenerate on the unstable
manifolds. Obviously, $\omega^k$ is also nondegenerate on the unstable
manifolds for all maps close to $T_A$.

It remains to show that $W^u$ carries a unique homology. We first
observe that the map $f$ is homotopic to the linear map $T_A$ and
hence the induced maps on the homology are exactly the same. By
Proposition \ref{prop30}, we have that the set of all homologies
carried by $W^u$ span an invariant subspace in $H_k(\T^n, \R)$. Since
the unstable manifold is expanded by approximately a factor of
$(\prod_{|\lambda_i| >1} \lambda_i)^n$.  Every eigenvectors of $f_*$
in this subspace has an eigenvalue close to $(\prod_{|\lambda_i| >1}
\lambda_i)^n$. However, and there is only one (up to a constant
multiple) eigenvector with the eigenvalue $(\prod_{|\lambda_i| >1}
\lambda_i)^n$. This implies that all the subsequential limit of the
currents is unique up to rescaling and the eigenvalue is exactly
$(\prod_{|\lambda_i| >1} \lambda_i)^n$.

\section{Proof of the main results}

In this section, we finish the proof of our main theorems.

First we define topological entropy using $(n, \ep)$-separated
sets. Let $f: M \ra M$ be a homeomorphism on a compact metric space
$M$. For any given positive integer $n$ and positive real number $\ep
>0$, a subset $S \subset M$ is said to be $(n, \ep)$-separated if for
any two distinct points $x, y \in S$, there is an integer $i$ with $0
\leq i \leq n$ such that $d_M(f^i(x), f^i(y)) \geq \ep$. Let $\#S$ be
the cardinality of the set $S$ and let
$$s(n, \ep) = \max\{ \#S \; | \; S \subset M \mbox{ is } (n,
\ep)\mbox{-separated} \}$$
we define
$$h(f, \ep) = \lim \sup_{n \ra \8} \frac{1}{n} \ln s(n, \ep).$$
and the topological entropy $$h_{top}(f)= h(f) = \lim_{\ep \ra 0^+}
h(f, \ep).$$

The topological entropy measures the growth of the trajectories. It is
certainly related to the geometric growth of the unstable foliations
in a partially hyperbolic system. Let $f \in \PH^r(M)$ be a partially
hyperbolic diffeomorphism and let $W^u$ be the unstable foliation of
$f$. Let $\chi_u(f) = \chi_{W^u}(f)$ be the geometric growth of $f$ on
the unstable foliation $W^u$. We have the following lemma.

\begin{lem}
$h(f) \ge {\chi_u(f)}$.
\end{lem}

\begin{proof}
For any given $\delta >0$, choose a point $x \in M$ and a small
$r >0$ such that $\xi_u(x, r) > \xi_u(f) - \delta/2$. By the
definition of $\chi_u(x, r)$, there exists
$N>0$ such that for all $n\geq N$ we have the following inequality
\[
\Vol(f^{n}(W^u_r(x)))>{e^{n(\chi_u(f)-\delta)}}
\]
For any $\ep >0$, we consider $(n,\epsilon)$-separated sets on the the
strongly unstable foliation $W^{u}(x)$.  Let $d_u(y, z)$ be the
distance between two points $y, z \in W^u(x)$ measured by the shortest
curve in the submanifold $W^u(x)$ between $y$ and $z$. Clearly $d_M(y,
z) \leq d_u(y, z)$, where $d_M$ is the distance between $y$ and $z$ on
$M$. Without loss of generality, we may assume that $\ep < r /(\max_{x \in
  M}\| Df_x \|)$. Let
$$\ep' = \inf\{ d_M(y, z) \; | \; y, z \in W^u(x); \; x \in M;\; \ep
\leq d_u(y, z) \leq r \}$$ We claim that $\ep' >0$. For otherwise, by
the continuity of the leaves, there exist a point $y \in M$ and a
sequence of points $z_i \in W^u(y)$ such that $d_M(y, z_i) \ra 0$, as
$i \ra \8$ and
$$\ep \leq d_u(y, z_i) \leq r, \mbox{ for all } i \in \N$$
This is impossible since the set $A= \{z \in W^u(y) \; \| \; \ep \leq
d_u(y, z) \leq r \}$ is compact and $d_M(y, z) >0$ for all $z \in
A$.

Let $S(n, \ep) \subset f^n(W^u_r(x))$ be a finite set such that
for any $x_i, x_j \in S(n, \ep)$, $x_i \neq x_j$, we have $d_u(y,
z) \geq \ep$. There is a constant $c_7>0$, depending on the
Riemannian metric and $k$, the dimension of the foliation $W^u$,
such that $\Vol (W^u_\ep(x_i)) \leq c_7 \ep^k$. Let $\#(S(n, \ep))$
be the cardinality of the set $S$. Then the total volume covered
by the $\ep$ balls around the points in $S(n, \ep)$ is less than
$\#(S(n, \ep))c_7\ep^k$. Since one can add points to $S(n, \ep)$ if
the total volume of these $\ep$ balls is less than the total
volume of $W^u_r(x)$, this implies that $S(n, \ep)$ can have at
least as many as points as $$\Vol(f^n(W^u_r(x)))/(c_7 \ep^k) $$

The pre-image $f^{-n}S(n, \ep) \subset W^u_r(x)$ is an $(n, \ep)$
separated set on the unstable foliation $W^u$. In fact, it is also an
$(n, \ep')$ separated set on $M$, where $\ep'$ is as defined
before. For given any two distinct points $y, z \in f^{-n}(S)$, we
have $d_u(f^n(y), f^n(z)) \geq \ep$ and $d_u(x, y) < r$. Since $\ep <
r /(\max_{x \in M}\| Df_x \|)$, there exists an integer $i$, $0 \leq i
\leq n$ such that $\ep \leq d_u(f^i(y), f^i(z)) \leq r$. Therefore,
$d_M(f^i(y), f^i(z)) \geq \ep'$.

Finally, \be h(f, \ep') &\geq& \lim \sup_{n \ra \8} \frac{1}{n} \ln
\#S(n, \ep) \nn \\
&\geq& \lim_{n \ra \8} \frac{1}{n} \ln \left(
e^{n(\chi_u(f) - \delta)}/(c \ep^k) \right) \nn \\
&=& \chi_u(f) - \delta \nn \ee

Since $\delta>0$ is arbitrary, $$h(f) = \lim_{\ep' \ra 0^+} h(f, \ep')
\geq \chi_u(f).$$
This completes the proof.
\end{proof}

Related to the topological entropy is the measure theoretic
entropy. Even though our results are about topological entropy, our
proof uses results from measure entropy, Lyapunov exponents and smooth
ergodic theory. Let $\nu$ be an invariant probability measure. Similar
to the definition of topological entropy, one can define an entropy,
$h_\nu(f) \geq 0$, associated with the invariant measure $\nu$, using
the so-called $(n, \ep)$-spanning set that covers a $\nu$ positive
measure set. We refer readers to Pollicott \cite{Po} and Robinson
\cite{Ro} for more details. However, later in the paper, we will use
the following equivalent definition.

One can define a {\it measure entropy} $h_{\nu}(f)$ as follows:
Call $\xi=\{A_1,...,A_r\}$ a (finite) {\it measurable partition}
of $X$ if the $A_i$ are disjoint measurable subsets of X covering
$X$. Now set
$$ H(\xi)=\overset{r}{\sum_{i=1}}\nu(A_i)\log\nu(A_i)$$
Then the limit \be h_\nu(f,\xi)&=&\lim_{n \ra \infty}\frac
{1}{n}H(\xi\vee f^{-1}\xi\vee\ldots\vee f^{-(n-1)}\xi) \nn \\
&=&\lim_{n \ra\infty}H(\overset{n-1}{\vee_{i=0}}f^{-i}\xi) \nn \ee
exists and one defines
$$ h_\nu(f)=\sup\{h_\nu(f,\xi):\ \xi\ \mbox{ is a finite measurable
partition of } X\}.$$ Let $\xi(A)=\{A_1,\ldots,A_k\}$ and
$\zeta(C)=\{C_1,\ldots,C_p\}$ be two finite partitions, we define the
{\it entropy of $\xi$ given to $\zeta$ }to be \be H(\xi | \zeta) & =
&-\overset{p}{\sum_{j=1}} \nu(C_j) \overset{k}{\sum_{i=1}}
\frac{\nu(A_i\cap
  C_j)}{\nu(C_j)}\log{\frac{\nu(A_i\cap C_j)}{\nu(C_j)}}\nn \\
&=&-\sum_{i,j}\nu(A_i\cap C_j)\log{\frac{\nu(A_i\cap C_j)}{\nu(C_j)}}
\nn \ee omitting the j-terms when $\nu(C_j)=0.$ Later in this paper we
shall use the following fact to compute $h_\nu(f,\xi)$:
$$ h_\nu(f,\xi)=\lim_{n \ra
\infty}H(\xi | (\overset{n}{\vee_{i=1}}f^{-i}\xi)).$$ (See
Walters~\cite{Wa1}, P82-83.)

From the definitions, it is easy to show that $h_\nu(f) \leq h(f)$.
Moreover, we have the following well known Theorem. We refer the
readers to Walters~\cite{Wa2} for excellent accounts.

\begin{thm} [Variational Principle] Let $M_{erg}$ be the set of all
  invariant ergodic measures, then $h(f)=\sup_{\nu\in
    M_{erg}}{{h_\nu(f)}}$. In other words, for all $\varepsilon>0$,
  there exists $\nu\in M_{erg}$ such that $ h_\nu(f)>h(f)-\varepsilon$.
\end{thm}

Another concept we will need to use is the Lyapunov exponents.  Let
$\nu$ be an invariant probability measure for $f \in \dfr$. For
$\nu$-a.e. $x \in M$, there exist real numbers $\lambda_1(x) > \ldots
> \lambda_l(x)$ ($l \leq n$); positive integers $n_1$, $\ldots$, $n_l$
such that $n_1 + \ldots + n_l =n$; and a measurable invariant
splitting $T_xM = E^1_x \oplus \cdots \oplus E^l_x$, with dimension
$\dim(E^i_x)=n_i$ such that $$\lim_{j \ra \8} \frac{1}{j} \log
\|D_xf^j(v_i) \| = \lambda_i(x),$$ whenever $v_i \in E_x^i$, $v \neq
0$.

These numbers $\lambda_1(x)$, $\cdots$, $\lambda_l(x)$ are called the
Lyapunov exponents for $x \in M$. If the probability measure $\nu$ is
ergodic, then these exponents are constants for $\nu$-a.e. $x \in
M$. The existence of these Lyapunov exponents is the result of
Oseledec's Multiplicative Ergodic Theorem.

An invariant measure $\nu$ is called {\it hyperbolic} on a invariant
set ${\Lambda}$ if $\nu({\Lambda}) > 0$ and $\nu$-a.e.~$x\in
{\Lambda}$ has the property that $\lambda_i(x) \neq 0$ for all $i=1,
\ldots, l$.

For a partially hyperbolic diffeomorphism, $f \in \PH(M)$, the
Lyapunov exponents can be relabeled in three groups, according where
their corresponding vectors are. We will write $\lambda^s_i$ for
Lyapunov exponents in $E^s$, $\lambda_i^c$ for exponents in $E^c$ and
$\lambda_i^u$ for exponents in $E^u$.

To study the continuity properties of the topological entropy, we now
consider diffeomorphisms close to a given $f \in \PH(M)$. Assume that
$W^u_f$, the unstable foliation for $f$, carries a unique non-trivial
homology, we say that $f$ {\it stably}\/ carries a unique non-trivial
homology if there is a neighborhood $V$ of $f$ in $\PH(M)$ such that
any $g \in V$, the unstable foliation $W^u_g$ uniquely carries the
same homology element (up to rescale). If the unstable foliation is
one dimensional and $f$ carries a unique non-trivial homology, it is
easy to show $f$ {\it stably}\/ carries a unique homology. We believe
that this is true in general, but we are not able to show
this. However, if all the subsequential limits of the currents $C_n$
are closed for $f$ and nearby maps, then one can show that $f$
carrying a unique nontrivial homology implies that $f$ stably carrying
a unique nontrivial homology. This is certainly true for maps on
$\T^n$ close to the linear one.

The homologies carried by the stable foliation is defined in the same
way by considering $f^{-1}$.

Under the assumption that $f$ stably carries a unique non-trivial
homology, the geometric expansion $\chi_u(f,x)$ is well-defined for
all $x\in M$ and is constant.  Furthermore $\chi_u(f)$ is locally
constant on $f$.

We now proceed with the proof of our main results. We will divide the
proof into several cases. Since the upper semi-continuity of
topological entropy is known from Yomdin's theorem for $C^\8$
diffeomorphisms, it suffices to show lower semi-continuity for our
results on dimension two.

In the proof we also need to consider $f^{-1}$. First we recall that
$h(f) = h(f^{-1})$. We can also define the volume
growth of the stable foliation $W^s$ under $f^{-1}$. Same as
$\chi_u(f)$, we can define $\chi_s(f) = \chi_u(f^{-1})$ and in the
same way, we have $h(f) \geq \chi_s(f)$.

\vspace{1ex}
\noindent {\bf Case 1}: Either $h(f) = \chi_u(f)$ or $h(f) = \chi_s(f)$.
\vspace{1ex}

This is a simple case. Assume $h(f) = \chi_u(f)$. By Proposition
~\ref{1}, $\chi_u(f)$ is locally constant, so there exists a
neighborhood of $V$ of $f$ in $\dfr$ such that for any $g \in V$,
we have $\chi_u(g)=\chi_u(f)$. Therefore
$h(g)\geq\chi_u(g)=\chi_u(f)=h(f)$. i.e., $h(f)$ is lower
semi-continuous. The case with $h(f) = \chi_s(f)$ is the same.

\vspace{1ex}
\noindent {\bf Case 2}: Both $h(f)>{\chi_u(f)}$ and $h(f)>{\chi_s(f)}$
hold.
\vspace{1ex}

This is our main case. Let $\nu$ be an ergodic invariant probability
measure for $f \in \dfr$ and let $\lambda_1, \ldots, \lambda_l$ be the
Lyapunov exponents associated with $\nu$. We have the following
Pesin-Ruelle inequality:
$$h_\nu(f) \leq \sum_{\lambda_i >0} \lambda_i.$$
For our purpose, we need a refined version of the Pesin-Ruelle
inequality, where we incorporate the geometric expansion $\chi_u(f)$
into the above formula. We have the following Theorem.

\begin{thm} \label{thmpr}
  Let $f \in \PH(M)$ be a partially hyperbolic diffeomorphism on a
  compact manifold $M$. Let $\nu$ be an ergodic measure and let
  $\lambda_i^c$ be the Lyapunov exponents corresponding to the center
  distribution $E^c$. Then the following estimate holds $${h_\nu(f)}
  \leq \sum_{\lambda_i^c>0}{\lambda_i^c} + {\chi_u(f)}.$$
\end{thm}

The proof of this theorem is quite involved. We postpone the proof to
the next section.

We return to the proof of the main theorem for the case with $h(f) >
\chi_u(f)$ and $h(f)>{\chi_s(f)}$. By
the variational principle, for any $\delta >0$ there is an ergodic
measure $\nu$ such that $$h_\nu(f) > h(f) - \delta.$$ Choosing
$\delta$ such that $$0< \delta \leq \min\{\frac{h(f) -
  \chi_u(f)}{3}, \frac{h(f) - \chi_s(f)}{3}\},$$
we have $$ h_\nu(f) > h(f) - \delta > \chi_u(f) + \delta.$$ By the
above proposition, for such measure $\nu$,
$$h_\nu(f) \leq \chi_u(f) + \sum_{\lambda_i^c>0}{\lambda_i^c}.$$
Therefore, $$\chi_u(f) + \sum_{\lambda_i^c>0}{\lambda_i^c} \geq
\chi_u(f) + \delta.$$ and therefore,
$$\sum_{\lambda_i^c>0}{\lambda_i^c} \geq
\delta > 0.$$ We can easily see that at least one of the
$\lambda_i^c$ must be larger than 0.

If $\dim E^c = 1$, then there is only one center exponent and
$\lambda^c >0$. Now consider $f^{-1}$ and recall that $h_\nu(f)=
h_\nu(f^{-1})$ and $\chi_s(f)=\chi_u(f^{-1})$, therefore
$h_\nu(f^{-1}) > \chi_u(f^{-1})$. Apply the same argument to
$f^{-1}$ and we get ${-\lambda^c}>0$. But that's a contradiction
to ${\lambda}>0$. Therefore we can't have both $h(f)
> \chi_u(f)$ and $h(f)>\chi_s(f)$. This is implies that we can
only have case 1 and $h(f)$ is actually the maximum of $\chi_u(f)$
and $\chi_s(f)$. But both these numbers are locally constant,
therefore $h(f)$ must be locally constant. This proves our theorem
for the case where $\dim E^c =1$.

Assume dim $E^c=2$. For the measure $\nu$, there are two center
Lyapunov exponents, $\lambda_1^c$ and $\lambda_2^c$. we may assume
that $\lambda_1^c \geq \lambda_2^c$. Above arguments show that
${\lambda_1^c}>0$. By considering $f^{-1}$, we have $- \lambda_2^c
>0$, or, $\lambda_2^c <0$. Since all other Lyapunov exponents are
nonzero, this implies that the measure $\nu$ is a hyperbolic ergodic
measure.

To complete our proof for $\dim E^c=2$, we need one more result
from Katok and Mendoza.

\begin{prop} [Katok and Mendoza] Assume that $\nu$ is an ergodic
  hyperbolic measure for a $C^{1+\alpha}$ diffeomorphism, $\alpha >0$.
  Then for any $\epsilon>0$, there exists a uniformly hyperbolic
  invariant set $\Lambda \subset M$ such that: $h(f|\Lambda)>
  {h_\nu(f)-\epsilon}$.
\end{prop}

The proof of this proposition can be found at Pesin(\cite{Pe},
P122-124).

Hyperbolic invariant sets persist under small perturbations. There is
a neighborhood $V$ of $f$ in $\dfr$ such that for any $g \in V$, there
is a hyperbolic invariant set $\Lambda_g$, close to $\Lambda$, such
that $g|_{\Lambda_g}: \Lambda_g \ra \Lambda_g$ is topologically
conjugate to $f|_\Lambda: \Lambda \ra \Lambda$. Therefore
$$h(g) \geq h(g|_{\Lambda_g}) = h(f|_\Lambda) > h_\nu(f) - \ep >
h(f)-\delta -\ep.$$ In other words, the entropy of $f$ is lower
semi-continuous.

This completes the proof of the main theorem, assuming Theorem
\ref{thmpr}. \qed

\section{A refined Pesin-Ruelle Formula}

In this section, we give a proof of Theorem \ref{thmpr}, a refined
Pesin-Ruelle formula for partially hyperbolic diffeomorphisms.

The usual approach to the proof relies on the fact that if a partition
of the manifold is fine enough, then the diffeomorphism, up to finite
number of iterates, can be approximated by its linearization and
therefore, the growth in the partition can be estimated by Lyapunov
exponents. However, the volume growth $\xi$ is the opposite of the
Lyapunov exponents, it gives estimates of volumes of large
surfaces. The difficulty is to incorporate these two, seemingly
opposite, concepts into the partitions.

Choose a small $ \delta > 0$,
for any given $\epsilon>0$, and any $x \in M$, there is an integer
$K_x$, depending on $x$, 
such that,
\be \Vol(f^i(W_{r}(x)))< \delta^k e^{i(\chi(f)+\epsilon)} \ee
for all $0< r \leq 10\delta $ and all $i \geq K_x$. For any positive
integer $K$, let $S_K$ be the set of points such that $K_x \leq
K$. Obviously, for any measure $\nu$ on $M$, we have $\nu(M\backslash
S_K) \ra 0$ as $K \ra \8$. We also observe that, there is a constant
$c_8$ such that for any $x \in M$, any positive integer $i$ and any $r
\leq \delta$,
\be \Vol(f^i(W_r(x)))< c_8 \delta^k (\sup_{x \in M} \|df\|)^{ki}.
\label{eqn30} \ee
Where $k$ is the dimension of the unstable foliation.

Fix a $m=lK$ where $l$ is a positive integer and let $B(y, t)$ be
a ball centered at $y$ with radius $t$. Since $M$ is
compact, there exists ${t_m}>0$ such that for every $0<t<{t_m}$,
$y \in M$, and $x\in{B(y,t)}$ we have
$$ {\frac{1}{2}{d_x}{f^m}({{\exp_{x}^{-1}{B(y,t))}} \subset
    {{\exp_{{f^m}x}}^{-1}f^m({B(y,t)})} \subset
    {{2d_x}{f^m}({\exp_x}^{-1}B(y,t))}}},$$
where $\exp_x$ is the exponential map at $x \in M$.

Now, for any chosen $\ep >0$, there is a positive number $\alpha >0$,
$\alpha < t_m /100$, such that for any partition $\xi$ with $\diam \xi
\leq 2\alpha$, we have:
$$h_\mu(f^m,\xi) \geq h_\mu(f^m)-\epsilon.$$

Let $d_u$ be the induced metric on $W(x)$, $x \in M$, from the
Riemannian structure on $M$. We introduce a dynamically defined new
metric on the manifold. Let $J$ be a positive integer such that the
following is true: for any point $x \in M$ and $y \in W(x)$, with
$d_u(x, y) \leq \delta$, we have $d_u(f^{-JK}(x), f^{-JK}(y)) \leq
\alpha$. Since $M$ is compact and the unstable leaves are uniformly
expanding, such integer $J$ exists. We now define a new metric
$d^J$ by
$$d^J(x, y) = d(f^{JK}(x), f^{JK}(y))  \alpha / \delta$$
and this metric also induces a metric, $d^J_u$, on the unstable
leaves, we have
$$d_u^J(x, y) = d_u(f^{JK}(x), f^{JK}(y))  \alpha / \delta.$$
An important property we have for this new metric is that $$d^J(x,
y) \leq  \alpha, \; \mbox{ whenever } \; d(f^{JK}(x), f^{JK}(y)) \leq 
\delta.$$

The metric $d^J$ depends on the choice of $J$, which is chosen to be a
large integer. A ball with metric $d^J$ is a thin tube-like
object. The center direction and the stable direction is very long and
the unstable direction is very short. On the unstable manifold, by
equation \eqref{eqn10}, we have
$$d_u^J(x, y) = d_u(f^{JK}(x), f^{JK}(y))  \alpha / \delta \geq
c_1^{-1} \lambda_u^{JK}d_u(x, y)  \alpha / \delta .$$
While on the center or center-stable manifold, if such manifolds do
exist,
$$d^J(x, y) \leq
c_1 (\lambda_c'')^{JK}d_u(x, y) \alpha / \delta .$$ This is also true on an
approximate center or center stable manifold with a slightly larger
$c_1$. Such approximate center and center stable manifold always
exist.  For a $d^J$ ball on $M$ with a small radius, the ratio of the
length in the center-stable direction to that of the unstable
direction is at most $2c_1^2(\lambda_c''/\lambda_u)^{JK}$.

For any fixed $\delta >0$ and small $\alpha >0$, we will choose $J$
such that \be c_1^{-1} \lambda_u^{JK} \alpha / \delta > 100, \; \mbox{
  and } \; c_1 (\lambda_c'')^{JK} \alpha / \delta < 1/100. \label{eqn45}
\ee We may increase $J$ by decreasing $\alpha$. Throughout this
section, $\alpha$ can be made arbitrarily small.

Finally we define a new metric $\rho$ on $M$ by
$$\rho(x, y) = d(x, y) + d_J(x, y),$$
for all $x, y \in M$. Clearly $\rho$ is a metric on $M$. Observe that,
by our choice of $J$ in equation \eqref{eqn45}, on the unstable
manifold, $d_u^J(x, y) \geq 100 d_u(x, y)$ and on an approximate center
or center stable manifold, $d^J(x, y) \leq d(x, y)/100$. Therefore we
have the following important property of the metric $\rho$: it is
dominated by the metric $d^J$ in the unstable direction and dominated
by $d$ in the center and stable direction. By equation \eqref{eqn45},
the set, we called it a $\rho$-ball, given by $$B_\rho (x, r) = \{
y \in M \; | \; \rho(x, y) \leq r \}$$ contains a $d^J_u$-ball in the
unstable direction with a radius $r/2$ and contains a regular, lower
dimensional ball of radius $r/2$ in the center-stable direction.

We continue our proof of Theorem \ref{thmpr}. The proof uses proper
partitions to estimate the entropy. There is a special partition of
the manifold $M$ which is described in the following statement.

\begin{lem} \label{lem25}
Given $\epsilon>0$, there is a partition $\xi$ of $M$ such that
\begin{enumerate}
\item $\diam \xi\leq 2 \alpha \leq {{t_m}/{50}}$ and therefore, ${h_\mu(f^m,
\xi)}\geq{h_\mu(f^m)-\epsilon}$;

\item for every element $C \in \xi$ there exist $\rho$-balls $B_\rho(x,r)$ and
$B_\rho(x,r^{\prime})$, such that $\alpha/4 < r' < r < \alpha$ and
$B_\rho(x,r^{\prime})\subset C \subset B_\rho(x,r)$;

\item there exists $0<r<{t_m}/20$ such that if $C \in \xi$ then
$C \subset B(y,r)$ for some $y\in M$, and if $x \in C$ then
\be \frac{1}{2} {d_x}{f^m} (\exp_x^{-1}B(y,r)) \subset \exp_{f^m
  x}^{-1} f^m C \subset 2d_x f^m (\exp_x^{-1} B(y,r)). \label{eqn55} \ee

\end{enumerate}
\end{lem}

\begin{proof}
  To construct such a partition, given $\alpha>0$, consider a maximal
  $2\alpha/3$-separated set $\Gamma$, with respect to the metric
  $\rho$. i.e., $\Gamma$ is a finite set of points for which ${\rho(x,
    y)}>2\alpha/3$ whenever $x$, $y\in\Gamma$, and for any point $z
  \in M$, there is a point $x \in \Gamma$ such that $\rho(x, z) \leq
  2\alpha/3$. For $x \in \Gamma$ set
\[
D_\Gamma(x) = \{ y\in{M} : \rho(y,x) \leq \rho(y,z), z \in \Gamma
\setminus\{x\} \}.
\] 
Obviously, $B_\rho(x,\alpha/3)\subset
{D_\Gamma(x)}\subset{B_\rho(x,2\alpha/3)}$, Note that the sets
$D_\Gamma(x)$ corresponding to different points $x\in\Gamma$ intersect
only along their boundaries, i.e., at a finite number of sub-manifolds
of co-dimension greater than zero. Since $\mu$ is a Borel measure, if
necessary, we can move the boundaries slightly so that they have zero
measure.

We may choose $\alpha$ arbitrarily small by increasing $J$. This
guarantees the properties in the lemma.
\end{proof}

\vspace{1ex}
Continuing with the proof of the theorem, observe that
\begin{equation}\label{four}
\begin{aligned}
h_\nu(f^m,\xi) =& \lim_{k\to\infty} H_\nu (\xi | f^m \xi \vee
\ldots \vee f^{km}\xi) \\
 \leq& H_\nu (\xi | f^m\xi) = \sum_{D\in f^m\xi} \nu(D) H(\xi|D)  \\
\leq& \sum_{D \in f^m\xi} \nu(D) \log \# \{ C\in\xi : C\cap D \neq
\emptyset\},
\end{aligned}
\end{equation}
where $H(\xi|D)$ is the entropy of $\xi$ with respect to conditional
measure on $D$ induced by $\nu$. To estimate the entropy, we need to
know the number of elements $C \in \xi$ that have nonempty
intersections with a given element $D \in f^m \xi$.  By Property 3 of
the partition $\xi$, we have a uniform control on the derivatives of
$f^m$ for each element $C \in \xi$. If we had used the regular metric
in our partition, then the estimate on the growth of the partition can
be easily obtained, this would leads to the standard Pesin-Ruelle
formula. In our case, we need to estimate the number of intersecting
partitions in two different directions: on the unstable direction and
on the center and stable directions.

We first estimate the growth of the partition in the unstable
direction. Consider the unstable disk $W_r(x)$ with $r \leq
\delta$. If $W_r(x) \cap S_K \neq \emptyset$, we have an estimate on
the $k$-dimensional volume of $f^K(W_r(x))$,
$$\Vol(f^K(W_r(x)))< \delta^k e^{K(\chi(f)+\epsilon)}.$$
Therefore, there is a constant $c_9$ such that $f^K(W_r(x))$ contains
at most $c_9e^{K(\chi(f)+\epsilon)}$ non-intersecting disks on the
unstable leave with radius not less than $\delta/20$. Similarly, for
any positive integer $i$,
$$\Vol(f^{iK}(W_r(x)))< \delta^k e^{iK(\chi(f)+\epsilon)}.$$
The set $f^{iK}(W_r(x))$ contains
at most $c_1e^{iK(\chi(f)+\epsilon)}$ non-intersecting disks on the
unstable leave with radius not less than $\delta/20$.

We now translate the above statements in terms of our new metric
$\rho$.  Let $D_x \subset W(x)$ be a piece of the unstable manifold
contained in a $d^J_u$-ball of radius $\alpha$, centered at
$x$. Suppose that $D_x \cap S_K \neq 0$. By the definition of $d^J_u$,
we have $f^{JK}(D_x) \subset W_\delta(f^{JK}(x))$. If $x \in S_K$,
then
$$\Vol(f^{(J+i)K} (D_x)) < \delta^k e^{iK(\chi(f)+\epsilon)}.$$
This implies that the $d^J_u$ volume, we denote it by $\Vol^J$, for
the set $f^K(D_x)$, is
$$\Vol^J(f^{K} (D_x)) = \Vol (f^{(J+1)K}(D_x)) < \delta^k
e^{K(\chi(f)+\epsilon)}.$$
and for any positive integer $i$,
$$\Vol^J(f^{iK} (D_x)) = \Vol (f^{(J+i)K}(D_x)) < \delta^k
e^{iK(\chi(f)+\epsilon)}.$$ Consequently, $f^{iK}(D_x)$ contains at
most $c_9e^{iK(\chi(f)+\epsilon)}$ non-intersecting disks on the
unstable leave with $d^J_u$ radius not less than $\delta/20$.

To summarize, if we partition the unstable leaves with sets which are
bounded between $d^J_u$ disks of radius $\alpha/20$ and $\alpha$, then
$f^{iK}$ image of any element of the partition covers at most
$c_1e^{iK(\chi(f)+\epsilon)}$ number of elements in that partition,
provided that the element contains a point in $S_K'$, where $S_k' =
f^{JK}(S_K)$.

Now we consider the partition $\xi$ in Lemma \ref{lem25}. By property
2 of the partition $\xi$, we have that for every element $C$ of $\xi$,
there is $x \in M$ such that $B_\rho(x, \alpha/4) \subset C \subset B_\rho(x,
\alpha)$. For any $D \in f^m\xi$, the $k$-dimensional growth in the
unstable direction is controlled the geometric growth. In the
center and stable directions, it is controlled by derivative of the
map $f$, since each element of the partition is bounded, from both
below and above, by balls with the normal metric in the center-stable
directions.
 
Before we proceed, we have the following simple lemma.

\begin{lem}\label{2}
There exists a constant $K_1>0$ such that for $D\in{f^m\xi}$,
$$\# \{ C\in\xi : C \cap D \neq \emptyset \} \leq K_1 \sup \{ \|
d_xf \|^{mn} : x \in M \},$$ where $n$ is the dimension
of the manifold.
\end{lem}

This can be shown by estimating the $\rho$ volume expansion of each
element in $C$ under $f^m$ and using property 2 of the
Lemma~\ref{lem25}. Each $C$ is bounded by a product of regular disks
in center-stable manifold (or approximate one) and a $\rho$ disk in
the unstable direction. On the unstable direction, the $\rho$ volume
expansion is bounded by, from inequality \eqref{eqn30}, $c_8 \delta^k
(\sup_{x \in M} \|df\|)^{km}$. The expansion of volume in center and
stable directions is bounded by the maximal of the derivative of the
map. The $\rho$ volume in the center-stable direction is close to the
real volume. Since the unstable foliation is absolutely continuous,
the lemma follows from Fubini theorem.\qed

\vspace{1ex}
We have a better exponential bound for the number of those sets $D$
such that $D=f(C') \in f^m\xi$ and $C'$ contain regular points for
the invariant measure $\nu$. More precisely, given $\epsilon>0$, let
$R_{m,\epsilon}$ be the set of forward regular points $x \in M$ which
satisfy the following condition: for $k>m$ and $v \in E_x^c$, \be
e^{k(\lambda(x,v)-\varepsilon)} \|v\| \leq \|d_xf^kv\| \leq
e^{k(\lambda(x,v)+\varepsilon)} \|v\|. \label{eqn65}\ee Here
$\lambda(x, v)$ is the Lyapunov exponent at $x$ corresponding to the
vector $v$,
$$\lambda(x, v) = \lim_{i \ra \8} \frac{1}{i} \ln \| d_xf^i v
\|.$$ The limit exists for $\nu$-a.e. $x \in M$.

Finally, we repartition every element of $\xi$ into two sets. For any
$C \in \xi$, let $$C^1 = \{ x \in C \;|\: W_\alpha(x)\cap S_K' \neq
\emptyset \}$$ and $C^2 = C\backslash C^1$. Let $\xi^1$ be the
collection of the sets of type $C^1$ and $\xi^2$ be the collection of
type $C^2$. Together $\xi^1$ and $\xi^2$ form a partition of the
manifold, we denote this new partition by $\xi'$.

The following lemma gives a estimate of number of intersections of
$f^m(\xi^1)$ with $\xi'$. The total measure for the sets in $\xi^2$ is
small and its contribution to the entropy will be given in another
estimate.

\begin{lem}\label{3}
  For any given $\epsilon>0$, there is a $N>0$, such that for any
  $m>N$, if $C^1 \in \xi^1$ and $D =f^m(C^1) \in f^m\xi'$ such that
  $C^1$ has a nonempty intersection with $R_m$, then there exists a
  constant $K_2>0$ such that
$$\# \{ C \in \xi' : C\cap{D} \neq \emptyset \} \leq K_2e^{\epsilon m}
e^{m(\chi_u(f) + \epsilon)}
\prod_{i:\lambda_i^c>0}e^{m(\lambda_i^c+\epsilon)}.$$
\end{lem}

\begin{proof}
To establish the inequality note that
\[
\#  \{C\in\xi' : C\cap{D} \neq \emptyset \} \leq 2\Vol_{\rho}
(B)(\diam_\rho \xi)^{-n},
\]
where $\Vol_\rho (B)$ denotes the $\rho$ volume of
\[
B = \{ y\in M : \rho (y, \exp_{f^m (x)}(d_xf^{m}(\exp_x^{-1} B^{\prime})))
< \diam_\rho \xi \}
\]
where $B' = B_\rho(x,2\diam_\rho C^{\prime})\cap S_K'$, $C' \in \xi^1$,
$f^{m}(C') = D$ and some $x \in C' \cap R_m$.  The set $B$ can be
thought as a fattened set $D$. Let $W^{cs}(x)$ be the center-stable
manifold of $x$. In fact, an approximate one will suffice. 
Let $E$ be the subset of $W^{cs}(x)$ such that
$$E = \{ y \in W^{cs}(x) \; | \; \rho(y, x) \leq 4\alpha \}.$$
Obviously, 
$$B' \subset \{ y \in M \; |\; \rho_u(y, z) \leq 3\diam_\rho
C^{\prime}, \; \mbox{ for some } z \in E\cap S_K' \}$$ i.e., $B'$ is
contained in the product of the set $E$ and unstable disks.  Since the
unstable foliation is absolutely continuous, by Fubini theorem, up to
a bounded factor, $\Vol_\rho (B')$ is bounded by the product of the
volume expansion of unstable disks and the volume expansion of $E$. By
the invariance of the unstable and center foliations, the same is true
for $\Vol_{\rho}(B)$. We have already obtained the $\rho$ volume
expansion in the unstable direction. For the set $E$, it is bounded by
a ball with the regular metric, whose tangent space is on
center-unstable direction (or arbitrarily close to the center-stable
direction). On both $E$ and $f^m(E)$, the metric $\rho$ is dominated
by the regular metric. $f^m(E)$ is approximately an
$(n-k)$-dimensional ellipsoid, whose total volume is bounded by the
product of the lengths of the axes. The length of the axis can be
estimated by $d_xf^{m}|_{E^c\oplus E^s}$, using equations
\eqref{eqn55} \eqref{eqn65}.  Those of the axes that correspond to
non-positive exponents are at most sub-exponentially larger. The
remaining axes are of size at most $e^{m(\lambda_i^c+\epsilon)}$, up
to a bounded factor, for all sufficiently large $m$. Therefore, \be
\Vol_\rho (B) &\leq& K_3 e^{m\epsilon} (\diam_\rho B)^{n}
e^{m(\chi_u(f) + \epsilon)}
\prod_{i : \lambda_i^c>0} e^{m(\lambda_i^c + \epsilon)} \nn \\
&\leq& K_3 e^{m\epsilon} (2\diam_\rho \xi)^{n}
e^{m(\chi_u(f)+\epsilon)} \prod_{i : \lambda_i^c>0} e^{m(\lambda_i^c +
  \epsilon)}, \nn \ee for some constant $K_3>0$. The lemma follows.
\end{proof}

By Lemmas ~\ref{2} and Lemma~\ref{3}, we obtain \be &&
mh_{\nu}(f)-\epsilon =
h_{\nu}(f^{m})-\epsilon \leq {h_{\nu}(f^{m},\xi)} \nn \\
&\leq& {\sum_{f^{-m}(D)=C' \in \xi^1, (C' \cap R_m) \neq \emptyset} \nu(D) (\log{K_2} +
  \epsilon{m}+{m\sum_{i : \lambda_i^c > 0}(\lambda_i^c + \epsilon)} +
  {m(\chi_u(f) + \epsilon)}}) \nn \\
&& + {\sum_{f^{-m}(D)=C' \in \xi^1, (C' \cap R_m) = \emptyset} \nu(D)
  (\log{2K_1} + nm \log
  \sup{\{\|d_x{f}\| : x\in{M}\}})} \nn \\
&& + \sum_{f^{-m}(D)=C' \in \xi^2} \nu(D)
  (\log{2K_1} + nm \log
  \sup{\{\|d_x{f}\| : x\in{M}\}})  \nn \\
&\leq & \log{K_2} + \epsilon{m} + m \sum_{i : \lambda_i^c>0}(\lambda_i^c +
\epsilon) + m(\chi_u(f) + \epsilon) \nn \\
&& + (\log{2K_1} + nm \log \sup{\{\|d_x{f}\| : x \in {M}\}})
\nu(M\setminus{ (R_m \cup S_K')}) \nn
\ee
By the multiplicative ergodic theorem, we have
$$\bigcup_{m\geq0}R_{m}(\epsilon) = M(mod0)$$ for every sufficiently
small $\epsilon > 0$. Since every point is in $S_K$ for some $K$ and
$\nu(S_K') = \nu(S_K)$, we
have $\nu (M\backslash S_K')
\ra 0$ as $K \ra \8$. It follows that
$$h_{\nu}(f) \leq {\epsilon + \sum_{i:\lambda_i^c>0}(\lambda_i^c +
  \epsilon)} + (\chi_u(f)+\epsilon).$$
Let $\epsilon \rightarrow 0$, we obtain the desired upper
bound.

This proves the theorem.

\section{Examples}

In this section, we give several examples where the main theorem fails
when we drop the assumption on homology.

Let $M_l$ be a compact orientable surface of genus $l \geq 2$, with
constant negative curvature. Let $g_t: SM_l \ra SM_l$ be the geodesic
flow on the unit tangent bundle of $M_l$. For any fixed $t >0$, $g_t$
is a partially hyperbolic diffeomorphism on $SM_l$. The stable,
unstable and center distributions are all one dimensional. The
topological entropy for $g_1$ is nonzero and for any $t \in \R$,
$h_{top}(g_t)= |t|h_{top}(g_1)$. Therefore in this case the topological
entropy is not locally constant. Obviously, the stable and unstable
foliations do not carry any nontrivial homology in this case. In fact,
any diffeomorphism that is isotopic to identity can not carry
non-trivial homology, since its induced action on homology is
trivial.

Now consider a map $f: SM_l \times S^1 \ra SM_l \times S^1$
defined in the following way. Let $\alpha: S^1=\R^1/\Z^1 \ra S^1$
be a diffeomorphism  which is close to identity and has three
fixed points $y_i=(i-1)/4$, $i=1, 2, 3$, which satisfy the
following:
$$\alpha^{\prime}(y_1)>1, \alpha^{\prime}(y_3)<1,
\alpha^{\prime}(y_2)= 1\ and\ \alpha^{\prime\prime}(y_2)\neq 0$$
Define $f(x,y)=(g_{1+\sin(2\pi y)}, \alpha(y))$. Then
\be h_{top}(f)&=&\max_ih_{top}(f|SM_g \times \{y_i\})\nn \\
&=&\max_{i}(1+\sin(2\pi y_i))h_{top}(g_1)=2h_{top}(g_1).\nn \ee
Now consider a family of diffeomorphisms
$$\alpha_{\epsilon}=\alpha(y,\epsilon)=\alpha+\epsilon,$$ and
$$f_\epsilon(x,y)=(g_{1+\sin(2\pi y)}, \alpha_\epsilon(y)).$$
For any fixed $\epsilon\geq 0$, $f_\epsilon$ is a partially
hyperbolic diffeomorphism on $SM_l\times S^1$ with $\dim E^u=\dim
E^s=1$ and $\dim E^c=2$.

Hyperbolic fixed points persist under small
perturbations, therefore for every $\epsilon$ sufficiently small,
$\alpha_\epsilon$ have unique fixed points $y_1^\epsilon$ and $
y_3^\epsilon$ close to $y_1$ and $y_3$ respectively.

For $i=2,$ since
$$ \alpha'(y_2)=1,
\alpha''(y_2)\neq 0, \and \ \frac{\partial
\alpha_\epsilon}{\partial\epsilon}(y_2,0)=1\neq 0,$$
saddle-node bifurcations occur at $(y_2,0)$. Therefore when
$\epsilon$ is small enough, on one side of $\alpha=
\alpha_0$(without loss of generality we can assume it's on the
left side), $\alpha_\epsilon$ has no fixed points close to $y_2$.
Therefore $\alpha_\epsilon$ has no other fixed points other than
$y_1^{\epsilon}$ and $y_3^{\epsilon}$, which implies

\be \lim_{ \epsilon \ra o^-}h_{top}(f_\epsilon)&=&\lim_{ \epsilon
\ra o^-}\max_{i=1,3} h_{top}(f|SM\times y_i^\epsilon)\nn\\
&=&h_{top}(g_1)\nn \ee

This means that the topological entropy of $f$ is not continuous.

%\bibliographystyle{plain}
%\bibliography{mybib}

\end{document}